\documentclass[11pt]{amsart}
\usepackage[T1]{fontenc}
\usepackage{lmodern}
\usepackage[margin=1.1in]{geometry}
\usepackage{amsmath,amssymb,amsthm,mathtools}
\usepackage{microtype}
\usepackage{xcolor}
\usepackage{aliascnt}
\usepackage[colorlinks=true,linkcolor=blue!50!black,citecolor=blue!50!black,urlcolor=blue!50!black]{hyperref}
\usepackage[capitalize,nameinlink]{cleveref}

\newcommand{\C}{\mathbb C}
\newcommand{\R}{\mathbb R}
\newcommand{\Z}{\mathbb Z}

\newcommand{\II}{\mathrm{II}_1}
\newcommand{\E}{E_N}                          
\newcommand{\Tr}{\operatorname{Tr}}           
\newcommand{\Proj}{\operatorname{Proj}}
\newcommand{\U}{\mathcal U}

\newcommand{\sa}{\mathrm{sa}}
\newcommand{\op}{\mathrm{op}}
\newcommand{\one}{\mathbf 1}                  
\newcommand{\Lam}{\Lambda}                    
\newcommand{\Palpha}{\mathcal P_\alpha(M)}    
\newcommand{\basic}[2]{\langle #1,e_{#2}\rangle}

\newcommand{\norm}[1]{\left\lVert#1\right\rVert}
\newcommand{\abs}[1]{\left\lvert#1\right\rvert}

\newif\ifshowtodos

\showtodosfalse  

\newcommand{\newsharedtheorem}[3][theorem]{%
  \newaliascnt{#2}{#1}\newtheorem{#2}[#2]{#3}\aliascntresetthe{#2}}
\newtheorem{theorem}{Theorem}[section]
\newsharedtheorem{lemma}{Lemma}
\newsharedtheorem{proposition}{Proposition}
\newsharedtheorem{corollary}{Corollary}
\newtheorem{mainthm}{Theorem}

\newsharedtheorem[mainthm]{maincor}{Corollary}
\theoremstyle{definition}
\newsharedtheorem{definition}{Definition}
\newsharedtheorem{example}{Example}
\theoremstyle{remark}
\newsharedtheorem{remark}{Remark}
\numberwithin{equation}{section}

\crefname{mainthm}{Theorem}{Theorems}
\crefname{maincor}{Corollary}{Corollaries}

\title[Relative dimension in extremal subfactors]%
{The relative dimension set of an extremal subfactor\\of finite index}
\author{Julio Cáceres}
\email{julio.e.caceres.gonzales@vanderbilt.edu}
\date{November 23, 2026}

\begin{document}

\begin{abstract}
Popa's relative dimension set $\Lam(M,N)$ of an inclusion of $\II$
factors $N\subset M$ consists of the scalars $\alpha$ for which some
projection $p\in M$ satisfies $E_N(p)=\alpha1$. Popa computed it for
index at most $4$, and for index $d>4$ outside the middle interval
$(t_d,1-t_d)$, where $t_d(1-t_d)=d^{-1}$; whether the middle interval
meets $\Lam(M,N)$ for every irreducible inclusion of index $d>4$ has
remained open. We show that for every extremal inclusion of finite
index $d>4$, the whole interval $[t_d,1-t_d]$ is contained in
$\Lam(M,N)$. Consequently, $\Lam(M,N)$ depends only on the index for
all extremal inclusions of finite index, in particular for all
irreducible ones. Extremality cannot be dropped: for $d>4$, Popa's
locally trivial inclusions are exactly those whose relative dimension
set misses the middle interval. Other non-extremal inclusions can
contain the whole middle interval in their relative dimension set, or
miss part of it.
\end{abstract}

\maketitle

\section{Introduction}\label{sec:intro}


In a $\II$ factor $M$, projections are classified up to equivalence by
their trace, and the trace has a module-theoretic meaning: $\tau(p)$ is
the $M$-dimension of the right $M$-module $pL^2(M)$. Every value in
$[0,1]$ occurs, so $M$ has a continuous dimension theory. For an
inclusion $N\subset M$ of $\II$ factors, $L^2(M)$ is also a right
$N$-module, and its $N$-dimension is the Jones index $[M:N]$
\cite{Jones83}. In contrast with dimensions of modules over a single
factor, the index is rigid: Jones proved that
\[
 [M:N]\in\{4\cos^2(\pi/n):n\ge3\}\cup[4,\infty]
\]
\cite{Jones83}.

Pimsner and Popa described the index through the conditional
expectation $\E\colon M\to N$ alone \cite{PP86} (see
\cite[1.2.2]{Popa94}): $[M:N]^{-1}$ is the
largest constant $\lambda\ge0$ with
\[
 \E(x)\ge\lambda x\qquad(x\in M_+).
\]
For a projection $p$ this says that $\E(p)\ge[M:N]^{-1}p$, so no
nonzero projection of $M$ looks smaller than $[M:N]^{-1}$ from the
point of view of $N$. In finite index the bound is sharp, and it is
attained by a projection whose expectation is a scalar
\cite[Prop.~1.5]{Popa89}: $M$ is then
the basic construction of a subfactor of $N$, implemented by a Jones
projection $e\in M$ with $\E(e)=[M:N]^{-1}1$; iterating this downward
basic construction gives the tunnel of $N\subset M$ \cite{Popa89}.

This led Popa to study the projections whose expectation onto $N$ is a
scalar \cite{Popa89}. Since $\tau$ is the expectation onto $\C1$, the
operator $\E(p)$ may be viewed as the dimension of $p$ relative to $N$,
and Popa's \emph{relative dimension set}
\[
 \Lam(M,N)=\{\alpha\in[0,1]:\ \E(p)=\alpha1\text{ for some }p\in\Proj(M)\}
\]
collects its scalar values. For $0<\alpha<1$, the condition
$\E(p)=\alpha1$ says that $\C p+\C(1-p)$ and $N$ form a commuting
square, which was Popa's original motivation.

Popa computed $\Lam(M,N)$ in terms of the index alone in most cases
\cite[Thms.~5.1--5.2]{Popa89}. If $[M:N]<4$, it is a finite set of
quotients of Jones polynomials evaluated at $[M:N]^{-1}$. If
$d=[M:N]\ge4$, let $t_d=(1-\sqrt{1-4/d})/2$ be the smaller root of
$t(1-t)=d^{-1}$, and define $a_0=0$ and $a_{k+1}=1/(d(1-a_k))$, an
increasing sequence with limit $t_d$. Then
\[
 \Lam(M,N)\cap[0,t_d]=\{a_k:k\ge0\}\cup\{t_d\},
\]
and symmetrically $\Lam(M,N)\cap[1-t_d,1]=\{1-a_k:k\ge0\}\cup\{1-t_d\}$.
These values come from towers of projections built from the tunnel.
What remained was the middle interval $(t_d,1-t_d)$, which is nonempty
when $d>4$. Popa asked whether $\Lam(M,N)$ meets it when $N'\cap M=\C1$,
and whether $\Lam(M,N)$ can be uncountable, is always closed, or can
contain an interval \cite[Problems~5.4.3 and 5.4.6]{Popa89}.

These questions have remained open. Popa returned to $\Lam(M,N)$ in
\cite[\S3]{Popa19}. He conjectured that $[t_d,1-t_d]\subseteq\Lam(M,N)$
for inclusions with principal graph $A_\infty$ and index $d>4$, noting
supporting evidence from the standard representation of the inclusion
\cite[\S3.3]{Popa19}; see also \cite{Popa23}. The set $\Lam(M,N)$ has
been related to angles
between intermediate subfactors \cite{BDLR19,BakshiGupta19}, and
rational points of the middle interval have been found for several
families of subfactors \cite{BakshiGuin25}. It remained unknown whether
$\Lam(M,N)$ is closed or uncountable, whether it can contain an
interval, or whether its middle part can contain an irrational number
\cite[p.~859]{BakshiGuin25}.

\subsection*{Main results}
Our main result fills in the middle interval for all extremal
inclusions. This completes Popa's computation for these inclusions, and proves his
conjecture from \cite[\S3.3]{Popa19} in a stronger form, since
inclusions with principal graph $A_\infty$ are irreducible, hence
extremal.

\begin{mainthm}\label{thm:A}
Let $N\subset M$ be an extremal inclusion of $\II$ factors with
$4<d=[M:N]<\infty$. Then
$[t_d,1-t_d]\subseteq\Lam(M,N)$. Consequently
\[
 \Lam(M,N)=\{a_k:k\ge0\}\cup[t_d,1-t_d]\cup\{1-a_k:k\ge0\}.
\]
\end{mainthm}

\begin{maincor}\label{cor:B}
Let $N\subset M$ be an inclusion of $\II$ factors with $d=[M:N]<\infty$.
\begin{enumerate}
\item If $d\le4$, then $\Lam(M,N)$ is the set computed by Popa
 \cite[Thms.~5.1--5.2]{Popa89} (\cref{thm:popa}).
\item If $d>4$ and $N\subset M$ is extremal, then $\Lam(M,N)$ is the set
 in \cref{thm:A}.
\end{enumerate}
In particular, $\Lam(M,N)$ depends only on the index for all extremal
inclusions of finite index, and hence for all irreducible inclusions of
finite index.
\end{maincor}

\Cref{cor:B} and its consequences in \cref{sec:final} answer several of
Popa's questions from \cite[\S5.4]{Popa89} for all extremal inclusions
of finite index, in particular for irreducible ones:
\begin{itemize}
\item the middle interval does not merely meet $\Lam(M,N)$
 (Problem~5.4.3); it is contained in it;
\item $\Lam(M,N)$ is always closed, and for $d>4$ it is uncountable and
 contains an interval (Problem~5.4.6);
\item $\Lam(pMp,pNp)=\Lam(M,N)$ for projections $p\in N$
 (Problem~5.4.1);
\item $\Lam(M_1,M)=\Lam(M,N)$, where $M_1$ is the basic construction
 (Problem~5.4.2).
\end{itemize}
These are collected in \cref{cor:popa-problems}. Popa also asked
whether such invariants depend only on the standard invariant
\cite[\S3.4]{Popa19}: for these inclusions, $\Lam(M,N)$ depends on the
index alone, while the analogous question for the sets of tuples
realized by partitions of unity remains open (\cref{rem:partitions}).

In particular, for every extremal inclusion with $d>4$ there is a
self-adjoint unitary $u\in M$ with $\E(u)=0$, which answers a question
Popa raised for inclusions with principal graph $A_\infty$
\cite[\S3.5]{Popa19} (\cref{rem:orthonormal-unitaries}).

Extremality cannot be dropped. An inclusion $N\subset M$ is \emph{locally
trivial} if $N=\{x+\theta(x):x\in fMf\}$ for a projection $f\in M$ with
$0<\tau(f)<1/2$ and an isomorphism $\theta\colon fMf\to(1-f)M(1-f)$
(\cref{def:locally-trivial}). These inclusions have index $d>4$, are not
extremal, and their relative dimension set misses the middle interval
\cite[Thm.~5.3]{Popa89} (\cref{rem:locally-trivial}). As Popa pointed
out, they are the only ones.

\begin{maincor}\label{cor:C}
Let $N\subset M$ be an inclusion of $\II$ factors with
$4<d=[M:N]<\infty$. Then $\Lam(M,N)\cap(t_d,1-t_d)=\emptyset$ if and only
if $N\subset M$ is locally trivial.
\end{maincor}

For inclusions that are neither extremal nor locally trivial, both
behaviors occur. Extremality is not necessary for the conclusion of
\cref{thm:A}: \cref{ex:tensor} is such an inclusion whose relative
dimension set is the set in \cref{thm:A}. On the other hand,
\cref{ex:gap} is such an inclusion that misses part of the middle
interval, so $\Lam(M,N)$ is not determined by the index for non-extremal
inclusions. It also has $\Lam(M_1,M)\ne\Lam(M,N)$, so the equality in
Problem~5.4.2 can fail without extremality (\cref{rem:problem-542}).

The restriction to finite index is natural. The set $\Lam(M,N)$ is the
scalar part of Popa's dilation problem \cite[\S7]{Popa89},
\cite[\S3.6]{Popa19}, which for maximal abelian subalgebras is the
Carpenter problem \cite{Kadison02a,Kadison02b,AK06}; in finite index the Pimsner--Popa
inequality rules out most positive contractions, so the scalar part is
the natural question. In a companion paper with Bisch we show that if
$N\subset M$ has uniform infinite index, for instance if
$N'\cap M=\C1$ and $[M:N]=\infty$, then every positive contraction of
$N$ is the expectation of a projection, so $\Lam(M,N)=[0,1]$; without
uniform infinite index, $\Lam(M,N)$ is not determined by the index
\cite{BCcarpenter}. Together with \cref{cor:B}, this shows that
$\Lam(M,N)$ is determined by $[M:N]$ for every irreducible inclusion of
$\II$ factors.

\subsection*{Strategy of the proof}
To prove \cref{thm:A}, fix $\alpha$ with $d\alpha(1-\alpha)>1$, and
consider the energy $F(p)=\norm{\E(p)-\alpha1}_2^2$ on the projections
of trace $\alpha$; then $\alpha\in\Lam(M,N)$ exactly when $F$ vanishes
somewhere. The idea is a dimension count. For matrices, a linear map
that is bounded below on $rM_k(\C)s$ cannot send it into a space of
smaller dimension. Suppose $p$ minimizes $F$ and put $a=\E(p)$. The
second variation of $F$ along rotations of $p$ shows that $\E$ is
bounded below on ``rectangles'' $rMs$, where $r\le p$ and $s\le1-p$ lie
in the high and the low spectral regions of $a$, respectively. These
rectangles are not modules over $M$ or $N$, so neither Hilbert-space
dimension nor the trace of the basic construction measures them. The
trace of Popa's symmetric enveloping algebra does: if $\E$ is bounded
below on $rMs$, then $\tau(r)\tau(s)\le d^{-1}\tau(e)\tau(f)$ whenever
$r\le e$ and $s\le f$ with $e,f\in\Proj(N)$ (\cref{lem:rectangle}). On
the other hand,
$\E(p)=a$ forces $p$ to carry at least an $\alpha$-fraction of the
upper tail of $a$, and $1-p$ at least a $(1-\alpha)$-fraction of its
lower tail (\cref{lem:balance}). Integrated over the spectrum, the two
estimates are compatible only if $\alpha(1-\alpha)\le d^{-1}$.

Two difficulties remain. First, $F$ need not attain its infimum. We
apply the Borwein--Preiss variational principle to $F^{1/4}$, which
gives a projection $p$ minimizing $F^{1/4}-2\tau(h\,\cdot\,)$ for some
$0\le h\le\eta$ (\cref{cor:BP}). With the exponent $1/4$, every error
term caused by the perturbation is a fixed small multiple of the main
terms, at every scale of $F$. So a single perturbed minimizer must
satisfy $F=0$. Working directly in $M$ avoids the usual detour through
the ultrapower $M^\omega$: there the infimum of $F$ is attained, but a
minimizer only shows that $\alpha$ is a limit of approximate solutions
in $M$. Second, the
perturbed minimizer commutes with $b=a-4F^{3/4}h$ rather than with $a$,
so the rectangles come from spectral projections of $b$, which need not
lie under projections of $N$. A spectral gap between the two spectral
regions lets us rotate them, by partial isometries close to the
identity, under spectral projections of $a$ (\cref{lem:sylvester}).
The endpoints $t_d$ and $1-t_d$ are covered by Popa's theorem.

\subsection*{Organization}
\Cref{sec:prelim} collects preliminaries: the basic construction and
extremality, Popa's computation of $\Lam(M,N)$, and the symmetric
enveloping trace. \Cref{sec:dimension} proves the dimension bound for
rectangles, \cref{sec:variations} constructs the perturbed minimizer
and computes its variations, and \cref{sec:proofA} combines them to
prove \cref{thm:A} and \cref{cor:B}. \Cref{sec:final} answers Popa's
questions and records some that remain open. \Cref{sec:nonextremal}
treats non-extremal inclusions: it proves \cref{cor:C} and constructs
\cref{ex:tensor,ex:gap}.

\subsection*{AI use statement}
Most of the proofs and arguments in this paper were developed with the
help of ChatGPT-6 Astra (OpenAI), which worked under the author's
guidance. The
writing and organization of the paper, including the checking of
references against the cited sources, were carried out with Claude
Opus 5.5 (Anthropic), with editorial direction from the author. The
author directed the project, made the mathematical and editorial
decisions, checked all proofs, and takes full responsibility for the
content of the paper.

\subsection*{Acknowledgements}
The author was partially supported by U.S. Army Research Office
(ARO) grant W911NF-23-1-0026 during this project. The author would also like to thank Dietmar Bisch for the many fruitful discussions that led to this problem, and Sorin Popa
for comments on a draft of this paper, in particular for the observation
of \cref{cor:C} and for suggesting the application to orthonormal
unitaries in \cref{rem:orthonormal-unitaries}.
\section{Preliminaries}\label{sec:prelim}

\subsection{Conventions}
All finite von Neumann algebras carry specified faithful normal tracial
states. For an inclusion $N\subset M$, the trace on $N$ is the restriction
of the trace $\tau$ on $M$, and $\E\colon M\to N$ is the trace-preserving
conditional expectation. We write $\norm{y}$ for the operator norm,
$\norm{y}_1=\tau(\abs y)$ and $\norm{y}_2=\tau(y^*y)^{1/2}$, and
$\norm{y}_{2,A}$ when the trace of an algebra $A$ must be specified;
$\Proj(M)$ is the set of projections of $M$. Inner products are linear
in the first variable; on $L^2(M)$,
$\langle\hat x,\hat y\rangle=\tau(y^*x)$. A corner $eMe$, for
$0\ne e\in\Proj(M)$,
carries the normalized trace $\tau_e(y)=\tau(y)/\tau(e)$.

For self-adjoint $y$, the positive and negative parts $y_+,y_-\ge0$
satisfy $y=y_+-y_-$ and $y_+y_-=0$. We write $y_+^2$ for $(y_+)^2$,
and similarly $\lambda_+^2=(\max(\lambda,0))^2$ for $\lambda\in\R$.

\subsection{The basic construction}
Let $N\subset M$ be an inclusion of $\II$ factors. On $L^2(M)$ let
$J\hat x=\widehat{x^*}$, and let $e_N$ be the orthogonal projection onto
$L^2(N)$. The basic construction is $M_1=\basic MN=(JNJ)'$. It carries a
unique normal semifinite faithful trace $\Tr$ with
$\Tr(xe_Ny)=\tau(xy)$ for $x,y\in M$. For a projection $T\in M_1$, the
range $TL^2(M)$ is a right $N$-module and
\[
 \Tr(T)=\dim_N\bigl(TL^2(M)\bigr),
\]
the right $N$-dimension; in particular $\Tr(1)=[M:N]$
\cite[Prop.~9.4.2, (9.3), Def.~8.5.4 and Def.~9.4.9]{AP}. Moreover
$M'\cap M_1=JMJ\cap(JNJ)'=J(N'\cap M)J$.

A finite family $m_1,\dots,m_k\in M$ is a \emph{Pimsner--Popa basis} if
$x=\sum_j m_j\E(m_j^*x)$ for all $x\in M$, equivalently
$\sum_j m_je_Nm_j^*=1$ as operators on $L^2(M)$. Such a family exists if
and only if $[M:N]<\infty$ \cite{PP86}; see also
\cite[Prop.~9.4.8]{AP}. In that case
\begin{equation}\label{eq:index-basis}
 [M:N]=\Tr(1)=\sum_j\Tr(m_je_Nm_j^*)=\sum_j\tau(m_jm_j^*).
\end{equation}
The index is invariant under compression by projections of the smaller
factor:
\begin{equation}\label{eq:jones-compression}
 [qMq:qNq]=[M:N]\qquad(0\ne q\in\Proj(N)),
\end{equation}
with normalized traces on both corners \cite{Jones83}; it also follows
from \cite[Exercise~8.9]{AP} and $\Tr|_M=[M:N]\,\tau$.

When $[M:N]<\infty$, the basic construction $M_1$ is a $\II$ factor, and so
is the commutant $N'=JM_1J$ of $N$ on $L^2(M)$
\cite[1.2.1(b)]{Popa94}; see also \cite[Prop.~9.4.2 and Def.~9.4.9]{AP}.
Let $\tau'$ denote the trace of $N'$.

\begin{definition}[\cite{PP86}]\label{def:extremal}
An inclusion $N\subset M$ of $\II$ factors with finite index is
\emph{extremal} if $\tau'(x)=\tau(x)$ for all $x\in N'\cap M$.
\end{definition}
For $0\ne r\in\Proj(N'\cap M)$, put $Nr=\{xr:x\in N\}$. The local index
formula
\begin{equation}\label{eq:local-index}
 [rMr:Nr]=[M:N]\,\tau(r)\,\tau'(r)
\end{equation}
holds \cite{Jones83}; see \cite[1.2.4(i)]{Popa94} and
\cite[Exercise~9.16]{AP}. Since $N'\cap M$ is finite-dimensional and
spanned by its projections, $N\subset M$ is extremal if and only if
$[rMr:Nr]=[M:N]\,\tau(r)^2$ for every nonzero $r\in\Proj(N'\cap M)$
\cite{PP86}; see \cite[1.2.5(i)--(ii)]{Popa94} and
\cite[\S2.2]{Popa23}. Pimsner and Popa's original formulation,
$E_{N'\cap M}(e)\in\C1$ for a Jones projection $e\in M$, is equivalent
\cite[Cor.~4.5]{PP86}. Irreducible inclusions ($N'\cap M=\C1$) are
extremal, and conversely an extremal inclusion with
$4<[M:N]<3+2\sqrt2$ is irreducible \cite{PP86}; see \cite[\S2.2]{Popa23}.

\subsection{The relative dimension set}
\begin{definition}[{\cite[Def.~1.1]{Popa89}}]\label{def:Lambda}
For an inclusion $N\subset M$ of $\II$ factors,
\[
 \Lam(M,N)=\{\alpha\in[0,1]:\ \E(p)=\alpha1\text{ for some }p\in\Proj(M)\}.
\]
\end{definition}
If $\E(p)=\alpha1$, trace preservation forces $\tau(p)=\alpha$, and
$\E(1-p)=(1-\alpha)1$. Thus $0,1\in\Lam(M,N)$ and $\Lam(M,N)$ is
invariant under $\alpha\mapsto1-\alpha$. For $0<\alpha<1$, the condition
$\E(p)=\alpha1$ says exactly that $\C p+\C(1-p)$ and $N$ form a
commuting square over $\C1$; this was Popa's original motivation
\cite{Popa89}. If $p\in\Proj(N'\cap M)$, then $\E(p)\in N'\cap N=\C1$,
so $\tau(\Proj(N'\cap M))\subseteq\Lam(M,N)$ \cite[Prop.~1.4]{Popa89}.
For irreducible inclusions this gives only $0$ and $1$, and every
projection realizing a value $0<\alpha<1$ fails to commute with $N$.

\begin{lemma}[Pimsner--Popa obstruction]\label{lem:PP}
If $[M:N]=d<\infty$, then $\Lam(M,N)\cap(0,d^{-1})=\emptyset$.
\end{lemma}
\begin{proof}
The Pimsner--Popa inequality gives $\E(a)\ge d^{-1}a$ for every $a\in M_+$
\cite{PP86}; see \cite[1.2.2]{Popa94}. If $\E(p)=\alpha1$ with
$0<\alpha<d^{-1}$, then $p\ne0$ by trace preservation, and compressing
the inequality by $p$ gives $\alpha p\ge d^{-1}p$, a contradiction.
\end{proof}
In fact $\inf(\Lam(M,N)\setminus\{0\})=[M:N]^{-1}$, with $\infty^{-1}=0$
\cite[Prop.~1.5]{Popa89}; the value $d^{-1}$ is attained by a Jones
projection of a downward basic construction \cite[1.2.2--1.2.3]{Popa94}.

For $d\ge4$ recall from the introduction
\begin{equation}\label{eq:an}
 t_d=\frac{1-\sqrt{1-4/d}}{2},\qquad
 a_0=0,\qquad a_{k+1}=\frac{1}{d(1-a_k)}\quad(k\ge0).
\end{equation}
Then $a_k$ increases to $t_d$, the smaller root of $t(1-t)=d^{-1}$.
The same recursion makes sense for $d=4\cos^2(\pi/(n+2))$, $n\ge1$, as
long as $a_k<1$; in that case $a_k<1$ for $k<n$ and $a_n=1$. Popa
writes these values as quotients of Jones polynomials evaluated at
$d^{-1}$ \cite[Thms.~5.1--5.2]{Popa89}; \eqref{eq:an} is the
corresponding recursion.

\begin{theorem}[Popa {\cite[Thms.~5.1--5.2]{Popa89}}]\label{thm:popa}
Let $N\subset M$ be an inclusion of $\II$ factors with $d=[M:N]<\infty$.
\begin{enumerate}
\item If $d=4\cos^2(\pi/(n+2))$ for some $n\ge1$, then
 $\Lam(M,N)=\{a_0,a_1,\dots,a_n\}$.
\item If $d\ge4$, then
 \[
 \Lam(M,N)\cap[0,t_d]=\{a_k:k\ge0\}\cup\{t_d\},\qquad
 \Lam(M,N)\cap[1-t_d,1]=\{1-a_k:k\ge0\}\cup\{1-t_d\}.
 \]
\end{enumerate}
No irreducibility is assumed.
\end{theorem}
By Jones' theorem \cite{Jones83}, an index $d<4$ has the form
$4\cos^2(\pi/(n+2))$ with $n\ge1$, so part~(1) applies. For $d=4$ the
intervals $[0,t_d]$ and $[1-t_d,1]$ in part~(2) cover $[0,1]$, since
$t_4=1/2$. Thus $\Lam(M,N)$ is determined by the index whenever
$d\le4$, with no hypothesis on the inclusion. For $d>4$ the
only undetermined part is $\Lam(M,N)\cap(t_d,1-t_d)$, which Popa asks
about in \cite[Problems~5.4.3 and 5.4.6]{Popa89}.

\subsection{The symmetric enveloping trace}
We regard $M$ as acting on $L^2(M)$ by left multiplication, and write
$y^{\op}=Jy^*J$ for right multiplication by $y\in M$, so that
$y^{\op}\hat x=\widehat{xy}$. Thus $M^{\op}=JMJ$ commutes with $M$, the
map $y\mapsto y^{\op}$ is linear with $(yz)^{\op}=z^{\op}y^{\op}$, and
$y^{\op}\ge0$ when $y\ge0$. Let
$\mathcal S_0=C^*(M,M^{\op},e_N)\subset B(L^2(M))$.

\begin{theorem}[Popa {\cite[Def.~1.8, Thm.~2.1]{Popa99}}]\label{thm:SEtrace}
Let $N\subset M$ be an extremal inclusion of $\II$ factors with finite
index. Then $\mathcal S_0$ has a tracial state.
\end{theorem}
Popa introduced the symmetric enveloping algebra in \cite{Popa94SE}.
His symmetric enveloping $C^*$-algebra is the universal $C^*$-algebra
generated by commuting copies of $M$ and $M^{\op}$ and a projection
implementing both $\E$ and $E_{N^{\op}}$, taken modulo the
representations satisfying a smoothness axiom \cite[Def.~1.8]{Popa99}.
The representation on $L^2(M)$ by left and right multiplication, with
$e_N$ the projection onto $L^2(N)$, satisfies this axiom, as noted
there, so $\mathcal S_0$ is a quotient of that algebra; and every such
quotient has a (unique) tracial state \cite[Thm.~2.1]{Popa99}; see also
\cite[Rem.~2.2.1]{Popa99}. The trace is built by averaging over unitaries
with the relative Dixmier property \cite[App.~A.1]{Popa99}, and its
construction requires neither amenability nor separability. The values
of the trace that we need follow from traciality alone
(\cref{lem:trace-formulas}). The trace need not be faithful on
$\mathcal S_0$; faithfulness is related to amenability of the standard
invariant \cite[\S5]{Popa99}, and it is never used below.

\section{A dimension bound for rectangles}\label{sec:dimension}

Throughout \cref{sec:dimension,sec:variations,sec:proofA}, $N\subset M$
is an extremal inclusion of $\II$ factors with finite index $d=[M:N]$,
and $\varphi$ is a tracial state on
$\mathcal S_0=C^*(M,M^{\op},e_N)\subset B(L^2(M))$, which exists by
\cref{thm:SEtrace}. Recall that $e_N$ commutes with $x$ and $x^{\op}$
for $x\in N$, and that $M$ and $M^{\op}$ commute.

Let $e,f\in\Proj(N)$ and $r,s\in\Proj(M)$ with $r\le e$ and $s\le f$.
\Cref{lem:rectangle} below compares two subspaces of $L^2(M)$: the
closure $rL^2(M)s$ of $rMs$, which is the range of $P=rs^{\op}$, and the
closure $eL^2(N)f$ of $eNf$, which is the range of $Q=ef^{\op}e_N$. Since
$\E(rxs)=e\,\E(rxs)\,f$, the expectation $\E$ maps the first space into
the second. Neither space is a module over $M$ or $N$: the first is
invariant only under $rMr$ on the left and $sMs$ on the right, the
second only under $eNe$ and $fNf$. Both are typically
infinite-dimensional. Hence neither Hilbert-space dimension nor the
dimension function $\Tr$ of the basic construction measures them. The
trace $\varphi$ does: $\varphi(P)=\tau(r)\tau(s)$ and
$\varphi(Q)=d^{-1}\tau(e)\tau(f)$ by \cref{lem:trace-formulas}.
\Cref{lem:rectangle} is then the analogue of the
fact that a linear map which is bounded below cannot send a space into
one of smaller dimension; for matrices, $rM_k(\C)s$ has dimension
$\operatorname{rank}r\cdot\operatorname{rank}s$.

\begin{lemma}[Trace formulas]\label{lem:trace-formulas}
Every tracial state $\varphi$ on $\mathcal S_0$ satisfies
\begin{align}
 \varphi(xy^{\op})&=\tau(x)\tau(y)\qquad(x,y\in M),\label{eq:trace-product}\\
 \varphi(e_N)&=d^{-1},\label{eq:trace-jones}\\
 \varphi(xy^{\op}e_N)&=d^{-1}\tau(x)\tau(y)\qquad(x,y\in N).\label{eq:trace-localized}
\end{align}
\end{lemma}
\begin{proof}
We use repeatedly that a positive functional $\psi$ on a $\II$ factor
$P$ with $\psi(uzu^*)=\psi(z)$ for all $u\in\U(P)$ and $z\in P$ is tracial
(unitaries span $P$), hence equal to $\psi(1)\tau_P$ by uniqueness of the
tracial state on $P$, which holds without any normality assumption
\cite[Prop.~4.1.3]{AP}. Note also that
$(u^{\op})^*y^{\op}u^{\op}=(uyu^*)^{\op}$ for $u\in\U(M)$.

\emph{Proof of \eqref{eq:trace-product}.} Fix $y\in M_+$. For $x\in M_+$
the operator $xy^{\op}$ is a product of commuting positive operators, so
$x\mapsto\varphi(xy^{\op})$ is a positive functional on $M$. For
$u\in\U(M)$ we have $(uxu^*)y^{\op}=u(xy^{\op})u^*$, so traciality of
$\varphi$ makes it unitarily invariant, and therefore
$\varphi(xy^{\op})=\varphi(y^{\op})\tau(x)$. Likewise $y\mapsto\varphi(y^{\op})$
is a positive functional on $M$, invariant under $y\mapsto uyu^*$, with
$\varphi(1)=1$; hence $\varphi(y^{\op})=\tau(y)$. Linearity gives
\eqref{eq:trace-product} for all $x,y\in M$.

\emph{Proof of \eqref{eq:trace-localized} up to the value of
$\varphi(e_N)$.} The same argument, using that $x$, $y^{\op}$ and $e_N$
commute for $x,y\in N$ and taking $u\in\U(N)$, shows
\begin{equation}\label{eq:trace-localized-c}
 \varphi(xy^{\op}e_N)=\varphi(e_N)\,\tau(x)\tau(y)\qquad(x,y\in N).
\end{equation}

\emph{Proof of \eqref{eq:trace-jones}.} Let $m_1,\dots,m_k$ be a
Pimsner--Popa basis, so that $\sum_jm_je_Nm_j^*=1$ on $L^2(M)$.
For $m\in M$, traciality, $e_N^2=e_N$ and $e_Nze_N=\E(z)e_N$ give
\[
 \varphi(me_Nm^*)=\varphi(e_Nm^*me_N)=\varphi(\E(m^*m)e_N)
 =\varphi(e_N)\,\tau(m^*m),
\]
the last step by \eqref{eq:trace-localized-c} with $y=1$. Summing over
$j$ and using $\tau(m_j^*m_j)=\tau(m_jm_j^*)$ and \eqref{eq:index-basis},
$1=\varphi(e_N)\sum_j\tau(m_j^*m_j)=d\,\varphi(e_N)$. Together with
\eqref{eq:trace-localized-c} this also proves \eqref{eq:trace-localized}.
\end{proof}

\begin{lemma}[Rectangular dimension comparison]\label{lem:rectangle}
Let $e,f\in N$ and $r,s\in M$ be projections with $r\le e$ and $s\le f$.
If, for some $\gamma>0$,
\begin{equation}\label{eq:bounded-below}
 \norm{\E(x)}_2\ge\gamma\norm{x}_2\qquad(x\in rMs),
\end{equation}
then
\begin{equation}\label{eq:rect-dimension}
 \tau(r)\tau(s)\le d^{-1}\tau(e)\tau(f).
\end{equation}
\end{lemma}
\begin{proof}
Put $P=rs^{\op}$ and $Q=ef^{\op}e_N$. These are projections in $\mathcal S_0$:
for $Q$, use that $e_N$, $e$ and $f^{\op}$ commute. Their ranges are
$rL^2(M)s$ and $eL^2(N)f$. Since $er=r$ and
$f^{\op}s^{\op}=(sf)^{\op}=s^{\op}$,
\[
 QP=e_N(er)(f^{\op}s^{\op})=e_NP .
\]
The vectors $\hat x$, $x\in rMs$, are dense in $PL^2(M)$, and
$e_N\hat x=\widehat{\E(x)}$. Hence \eqref{eq:bounded-below} gives
$Pe_NP\ge\gamma^2P$ in $B(L^2(M))$, and so in the $C^*$-algebra
$\mathcal S_0$.

Consequently $Pe_NP$ is invertible in the corner $P\mathcal S_0P$; let
$T\in P\mathcal S_0P$ be its inverse there, and put
$V=e_NPT^{1/2}\in\mathcal S_0$. Then
\[
 V^*V=T^{1/2}(Pe_NP)T^{1/2}=P,
\]
so $V$ is a partial isometry, and $V=QPT^{1/2}=QV$ gives
$VV^*=QVV^*Q\le Q$. Since $\varphi$ is a tracial state,
\[
 \varphi(P)=\varphi(V^*V)=\varphi(VV^*)\le\varphi(Q).
\]
By \cref{lem:trace-formulas}, $\varphi(P)=\tau(r)\tau(s)$ and
$\varphi(Q)=d^{-1}\tau(e)\tau(f)$.
\end{proof}

\section{A perturbed minimizer and its variations}\label{sec:variations}

\subsection{The energy}
Fix $0<\alpha<1$ and let
\[
 \Palpha=\{p\in\Proj(M):\tau(p)=\alpha\},\qquad
 F(p)=\norm{\E(p)-\alpha1}_2^2 .
\]
Since $\E(p)=\alpha1$ forces $\tau(p)=\alpha$, we have
$\alpha\in\Lam(M,N)$ if and only if $F$ vanishes somewhere on $\Palpha$.
For $p\in\Palpha$, trace preservation and $0\le\E(p)\le1$ give
\begin{equation}\label{eq:F-bound}
 F(p)=\tau(\E(p)^2)-\alpha^2\le\tau(\E(p))-\alpha^2=\alpha(1-\alpha)\le\tfrac14 .
\end{equation}
We regard $\Palpha$ as a metric space with the metric
$\operatorname{dist}(r,s)=\norm{r-s}_2$. On $\Palpha$,
\begin{equation}\label{eq:distance-overlap}
 \norm{r-s}_2^2=2\alpha-2\tau(rs).
\end{equation}

\begin{lemma}\label{lem:complete}
$(\Palpha,\operatorname{dist})$ is a complete metric space, and $F$ is
continuous on it.
\end{lemma}
\begin{proof}
The operator-norm unit ball of $M$ is complete for $\norm{\cdot}_2$
\cite[Prop.~2.6.4]{AP}, so a Cauchy sequence $(p_n)$ in $\Palpha$ converges in
$\norm{\cdot}_2$ to some $p\in M$ with $\norm p\le1$. Then $p=p^*$, since
the adjoint is $\norm{\cdot}_2$-isometric, and
$\norm{p_n^2-p^2}_2\le\norm{p_n(p_n-p)}_2+\norm{(p_n-p)p}_2\le2\norm{p_n-p}_2$
gives $p^2=p$; also $\tau(p)=\lim\tau(p_n)=\alpha$. Continuity of $F$
follows from $L^2$-contractivity of $\E$.
\end{proof}

\subsection{A perturbed minimizer}
The function $F$ need not attain its infimum on $\Palpha$, and we do not
know a priori that the infimum is zero. A variational principle
produces an exact minimizer after an arbitrarily small perturbation. We
apply it to $F^{1/4}$ rather than to $F$: with this exponent the
perturbation's effect on the first variation is of order $F^{3/4}$, which
keeps every error term in \cref{sec:proofA} proportional to the main
terms.

Recall that a \emph{gauge-type function} on a metric space
$(X,\operatorname{dist})$ is a continuous
$\rho\colon X\times X\to[0,\infty)$ with $\rho(x,x)=0$ such that for
every $\epsilon>0$ there is $\epsilon'>0$ with $\rho(y,z)\le\epsilon'$
implying $\operatorname{dist}(y,z)<\epsilon$.

\begin{theorem}[Borwein--Preiss principle, metric form]\label{thm:BP}
Let $(X,\operatorname{dist})$ be a complete metric space,
$f\colon X\to\R$ lower semicontinuous and bounded below, $\rho$ a
gauge-type function on $X$, and $\mu_i>0$ for $i\ge0$. Let $\epsilon>0$
and let $x_0\in X$ satisfy $f(x_0)\le\inf_Xf+\epsilon$. Then there are
$\bar x\in X$ and $x_1,x_2,\ldots\in X$ with $x_i\to\bar x$ such that
$\sum_{i\ge0}\mu_i\rho(\bar x,x_i)<\infty$ and
\[
 f(x)+\sum_{i\ge0}\mu_i\rho(x,x_i)\ \ge\ f(\bar x)+\sum_{i\ge0}\mu_i\rho(\bar x,x_i)
 \qquad(x\in X).
\]
\end{theorem}
This is the metric version \cite{LiShi00} of the Borwein--Preiss smooth
variational principle \cite{BorweinPreiss87}; see also
\cite[Thm.~2.5.2]{BorweinZhu05}. In the stated form it is the
infinite-series case of the extended principle of \cite{KPS}, in the
corrected form
\cite[Thm.~3 and Rem.~5(1)]{KPScorr}, with any sequence
$\epsilon_i\downarrow0$.

\begin{corollary}[Perturbed minimizer]\label{cor:BP}
For every $\eta>0$ there are $p\in\Palpha$ and $h\in M$ with
$0\le h\le\eta1$ such that $p$ minimizes
\begin{equation}\label{eq:functional}
 J(r)=F(r)^{1/4}-2\tau(hr)\qquad(r\in\Palpha).
\end{equation}
\end{corollary}
\begin{proof}
Apply \cref{thm:BP} to $X=\Palpha$ with the metric $\norm{r-s}_2$, the
continuous function $f=F^{1/4}$, the gauge-type function
$\rho(r,s)=\norm{r-s}_2^2$, $\mu_i=\eta\,2^{-i-1}$, $\epsilon=1$ and any
$x_0\in\Palpha$; by \eqref{eq:F-bound}, $0\le f\le1$, so
$f(x_0)\le\inf f+\epsilon$. Let $p=\bar x$ and $p_i=x_i$ ($p_0=x_0$). By
\eqref{eq:distance-overlap},
\[
 \sum_{i\ge0}\mu_i\norm{r-p_i}_2^2=2\alpha\eta-2\tau(hr),\qquad
 h=\sum_{i\ge0}\mu_ip_i ,
\]
where the series for $h$ converges in norm and $0\le h\le\sum_i\mu_i=\eta$.
Thus $p$ minimizes $F^{1/4}-2\tau(h\,\cdot\,)$ on $\Palpha$.
\end{proof}

\subsection{First and second variation}
\begin{lemma}[Stationarity and the rectangular Hessian]\label{lem:variations}
Let $p\in\Palpha$ minimize \eqref{eq:functional} with $0\le h\le\eta1$,
and suppose $F=F(p)>0$. Put
\[
 a=\E(p),\qquad q=1-p,\qquad b=a-4F^{3/4}h .
\]
Then
\begin{equation}\label{eq:commute}
 [p,b]=0,\qquad \norm{a-b}\le4\eta F^{3/4},
\end{equation}
and
\begin{equation}\label{eq:hessian}
 \norm{\E(x)}_2^2\ge\tau(bxx^*)-\tau(bx^*x)\qquad(x\in pMq).
\end{equation}
\end{lemma}
\begin{proof}
Let $K\in M$ with $K^*=-K$ and $p(t)=e^{tK}pe^{-tK}\in\Palpha$; dots denote
derivatives at $t=0$, so $\dot p=[K,p]$ and $\ddot p=[K,[K,p]]$. Write
$F(t)=F(p(t))=\tau(a(t)^2)-\alpha^2$ with $a(t)=\E(p(t))$. Since $a\in N$,
$\tau(a\E(y))=\tau(ay)$ for $y\in M$, and so
\[
 \dot F=2\tau(a\dot p),\qquad \ddot F=2\norm{\E(\dot p)}_2^2+2\tau(a\ddot p).
\]
As $F>0$, the function $t\mapsto J(p(t))=F(t)^{1/4}-2\tau(hp(t))$ is
smooth near $0$ and has a minimum at $t=0$.

\emph{First variation.} With $b=a-4F^{3/4}h$,
\[
 0=\tfrac{d}{dt}J(p(t))\big|_{t=0}
 =\tfrac14F^{-3/4}\dot F-2\tau(h\dot p)
 =\tfrac12F^{-3/4}\tau(b[K,p])=\tfrac12F^{-3/4}\tau([p,b]K).
\]
The commutator $C=[p,b]$ satisfies $C^*=-C$; taking $K=C^*$ gives
$\norm C_2^2=0$, so $[p,b]=0$. The norm bound follows from
$0\le h\le\eta1$.

\emph{Second variation.} At the minimum,
\[
 0\le\tfrac{d^2}{dt^2}J(p(t))\big|_{t=0}
 =\tfrac14F^{-3/4}\ddot F-\tfrac3{16}F^{-7/4}\dot F^2-2\tau(h\ddot p).
\]
Multiplying by $2F^{3/4}$ and inserting the formula for $\ddot F$,
\[
 \norm{\E(\dot p)}_2^2+\tau(b\ddot p)\ge\frac{3}{8F}\dot F^2\ge0 .
\]
Given $x\in pMq$, take $K=x^*-x$. Using $px=x$, $xp=0$ and $x^2=0$, one
finds $\dot p=x+x^*$ and $\ddot p=2(x^*x-xx^*)$, so
\[
 \norm{\E(x+x^*)}_2^2\ge2\bigl(\tau(bxx^*)-\tau(bx^*x)\bigr).
\]
The right side is unchanged if $x$ is replaced by $ix$. With $y=\E(x)$,
adding the two inequalities and using
$\norm{y+y^*}_2^2+\norm{y-y^*}_2^2=4\norm{y}_2^2$ proves
\eqref{eq:hessian}.
\end{proof}

\begin{remark}
The rotation generated by $x^*-x$, with $x\in pMq$, changes $\E(p)$ to
first order by $\E(x)+\E(x)^*$, and the rotation generated by
$(ix)^*-ix$ changes it by $i(\E(x)-\E(x)^*)$; these are twice the real
part and minus twice the imaginary part of $\E(x)$, which is in general
not self-adjoint. The gain $\tau(bxx^*)-\tau(bx^*x)$ is the same for
both rotations, so adding the two second-order inequalities bounds
$\norm{\E(x)}_2^2$, the quantity to which \cref{lem:rectangle} applies.
A single rotation gives the same bound up to a factor $2$, which would
only change the constants in \eqref{eq:bufferconditions}. The
term $\tau(bxx^*)-\tau(bx^*x)$ is the second-order gain from moving mass
of $p$ from high to low spectral regions of $b$. Inequality
\eqref{eq:hessian} says that at a minimum this gain is paid for by the
first-order visible change.
\end{remark}

\section{Proof of Theorem A}\label{sec:proofA}

Let $p$ be the minimizer from \cref{cor:BP} and suppose $F=F(p)>0$. By
\cref{lem:variations}, $p$ commutes with $b=a-4F^{3/4}h$, where
$a=\E(p)$, the distance $\kappa:=\norm{a-b}$ is at most $4\eta F^{3/4}$,
and the Hessian inequality \eqref{eq:hessian} holds.

We derive a contradiction by comparing an index side and a density side,
as in the matrix heuristic of \cref{sec:dimension}; $2$-norms of
spectral tails play the role of dimensions. Write $q=1-p$.
\begin{itemize}
\item \emph{Index side} (\cref{lem:tails}). For a \emph{threshold}
 $c\in\R$, a real number at which we split the spectrum of $a$, and a
 \emph{buffer} $\ell>0$,
 \[
 \norm{(b-c-\ell)_+p}_2\,\norm{(c-\ell-b)_+q}_2
 \le d^{-1/2}\,\norm{(a-c)_+}_2\,\norm{(c-a)_+}_2 .
 \]
 Inequality \eqref{eq:hessian} bounds $\E$ below on rectangles $rMs$
 with $r=p\,\one_{(c+\ell+t,\infty)}(b)$ and
 $s=q\,\one_{(-\infty,c-\ell-t')}(b)$; these are projections because $p$
 commutes with $b$. \Cref{lem:rectangle}, however, needs $r$ and $s$ to
 lie under projections of $N$, such as the spectral projections
 $\one_{(c+t,\infty)}(a)$ and $\one_{(-\infty,c-t')}(a)$ of $a$, and since
 $b\ne a$ this fails in general. We therefore rotate $r$ and $s$, by
 partial isometries close to the identity, to projections lying under
 these spectral projections of $a$. The rotation preserves traces and
 nearly preserves the lower bound for $\E$. The buffer $\ell$ provides
 the spectral gap that makes the rotation small (\cref{lem:sylvester}).
 If $h=0$, then $b=a$ and no rotation is needed.
\item \emph{Density side} (\cref{lem:balance}). Since $\E(p)=a$,
 \[
 \norm{(a-c)_+p}_2^2\ge\alpha\norm{(a-c)_+}_2^2,\qquad
 \norm{(c-a)_+q}_2^2\ge(1-\alpha)\norm{(c-a)_+}_2^2 ,
 \]
 and for a balanced threshold $c$ both norms on the right equal
 $\sqrt H$ with $H\ge F/2$.
\end{itemize}
The product of the densities, $\alpha(1-\alpha)$, then exceeds the
index bound $d^{-1}$, up to errors of size $\kappa+\ell$. With
$\ell=\delta\sqrt F$ and $\eta=\delta^2/16$ every error is a fixed small
multiple of the main terms, at every scale $F>0$.

\subsection{Buffered tails}
The tails of $b$ will be rotated under spectral projections of $a$. The
size of the required rotation is controlled by the following estimate,
a form of the Davis--Kahan $\sin\Theta$ theorem \cite{DavisKahan70}.

\begin{lemma}[Sylvester bound]\label{lem:sylvester}
Let $g,r\in\Proj(M)$, let $A\in gMg$ and $B\in rMr$ be self-adjoint with
$A\le\lambda g$ and $B\ge(\lambda+\ell)r$ for some $\lambda\in\R$ and
$\ell>0$, and let
$X,Y\in gMr$ satisfy $AX-XB=Y$. Then $\lVert X\rVert\le\lVert Y\rVert/\ell$.
\end{lemma}
\begin{proof}
For $v\ge0$ let $X(v)=e^{vA}Xe^{-vB}$, with the exponentials computed in
the corners $gMg$ and $rMr$, so that $X(0)=X$. Since $A$ commutes with
$e^{vA}$ and $B$ with $e^{-vB}$, the product rule gives
$X'(v)=e^{vA}(AX-XB)e^{-vB}=e^{vA}Ye^{-vB}$. The spectral bounds give
$\lVert e^{vA}\rVert\,\lVert e^{-vB}\rVert\le e^{-v\ell}$, so
$X(V)\to0$ as $V\to\infty$, and letting $V\to\infty$ in
$X(V)-X=\int_0^VX'(v)\,dv$ yields
\[
 X=-\int_0^\infty e^{vA}Ye^{-vB}\,dv,\qquad
 \lVert X\rVert\le\lVert Y\rVert\int_0^\infty e^{-v\ell}\,dv=\frac{\lVert Y\rVert}{\ell}.
 \qedhere
\]
\end{proof}

With the conventions of \cref{sec:prelim} for positive parts, we use the
layer-cake identity
\begin{equation}\label{eq:layer-cake}
 \lambda^2=\int_0^\infty2t\,\one_{(t,\infty)}(\lambda)\,dt\qquad(\lambda\ge0),
\end{equation}
applied through spectral calculus.

\begin{lemma}[Buffered-tail estimate]\label{lem:tails}
Let $p$ be a projection, $q=1-p$, $a=\E(p)$, and $b\in M$ self-adjoint
with $[p,b]=0$, $\norm{a-b}\le\kappa$, and satisfying \eqref{eq:hessian}
for all $x\in pMq$.
Let $\ell>0$ with
\begin{equation}\label{eq:bufferconditions}
 \theta:=\frac{\kappa}{\ell}\le\frac12,\qquad 4\theta<\sqrt{2\ell}.
\end{equation}
Then for every $c\in\R$,
\begin{equation}\label{eq:tails}
 \norm{(b-c-\ell)_+p}_2\,\norm{(c-\ell-b)_+q}_2
 \le d^{-1/2}\,\norm{(a-c)_+}_2\,\norm{(c-a)_+}_2 .
\end{equation}
\end{lemma}
\begin{proof}
We rotate the tails of $b$ under spectral projections of $a$ (Steps~1
and~2), check that the lower bound for $\E$ survives the rotation
(Step~3), apply \cref{lem:rectangle} (Step~4), and integrate (Step~5).

For $t\ge0$ define
\begin{align*}
 r_t&=p\,\one_{(c+\ell+t,\infty)}(b),&
 e_t&=\one_{(c+t,\infty)}(a),\\
 s_t&=q\,\one_{(-\infty,c-\ell-t)}(b),&
 f_t&=\one_{(-\infty,c-t)}(a).
\end{align*}
Since $p$ commutes with $b$, the $r_t$ and $s_t$ are projections
commuting with $b$; and $e_t,f_t\in N$. The upper and lower tails will
be paired independently: we compare $r_t$ with $s_{t'}$ for all
$t,t'\ge0$, since the final integration runs over the whole quadrant.
The tails of $b$ are cut a distance $\ell$ beyond the corresponding
tails of $a$.

\emph{Step 1: rotating the tails under projections of $N$; angles.} We show $\lVert(1-e_t)r_t\rVert\le\theta$. Write
$e=e_t$, $r=r_t$, $\lambda=c+t$. In the corner algebras $(1-e)M(1-e)$ and $rMr$,
the elements $A=(1-e)a(1-e)$ and $B=rbr$ satisfy $A\le\lambda(1-e)$, because
$1-e=\one_{(-\infty,\lambda]}(a)$, and $B\ge(\lambda+\ell)r$: a spectral gap of
$\ell$. The element $X=(1-e)r$ satisfies the Sylvester equation
$AX-XB=Y$ with $Y=(1-e)(a-b)r$ and $\lVert Y\rVert\le\kappa$, since $1-e$
commutes with $a$ and $r$ with $b$. By \cref{lem:sylvester} (with
$g=1-e$), $\lVert X\rVert\le\lVert Y\rVert/\ell\le\theta$. This argument
uses only that $r_t$ is a projection commuting with $b$ on which
$b\ge c+\ell+t$. Running it with $-a,-b,-c$ in place of $a,b,c$ and $q$
in place of $p$, which turns $r_t$ into $s_t$ and $e_t$ into $f_t$, gives
$\lVert(1-f_t)s_t\rVert\le\theta$ for all $t\ge0$.

Now fix $t,t'\ge0$.

\emph{Step 2: rotating the tails under projections of $N$; partial
isometries.} If $r_t\ne0$, then
$r_te_tr_t=r_t-r_t(1-e_t)r_t\ge(1-\theta^2)r_t$, so $r_te_tr_t$ is
invertible in $r_tMr_t$. The polar decomposition
$e_tr_t=w\,(r_te_tr_t)^{1/2}$ gives a partial isometry $w$ with
$w^*w=r_t$ and $ww^*\le e_t$, and
\[
 \lVert w-r_t\rVert\le\lVert w\bigl(r_t-(r_te_tr_t)^{1/2}\bigr)\rVert
 +\lVert e_tr_t-r_t\rVert\le\bigl(1-\sqrt{1-\theta^2}\bigr)+\theta\le2\theta .
\]
Likewise there is $v$ with $v^*v=s_{t'}$, $vv^*\le f_{t'}$ and
$\lVert v-s_{t'}\rVert\le2\theta$. (If $r_t$ or $s_{t'}$ is zero, the
inequality \eqref{eq:distribution} below is trivial.)

\emph{Step 3: the lower bound survives the rotation.} For $y\in r_tMs_{t'}
\subset pMq$, the spectral cuts give $\tau(byy^*)\ge(c+\ell+t)\norm y_2^2$
and $\tau(by^*y)\le(c-\ell-t')\norm y_2^2$, so by \eqref{eq:hessian}
\begin{equation}\label{eq:tail-lower}
 \norm{\E(y)}_2^2\ge(2\ell+t+t')\norm y_2^2 .
\end{equation}
Put $\tilde r_t=ww^*\le e_t$ and $\tilde s_{t'}=vv^*\le f_{t'}$. Every
$x\in\tilde r_tM\tilde s_{t'}$ equals $wyv^*$ with $y=w^*xv\in r_tMs_{t'}$, and
$\norm x_2=\norm y_2$. Since $x-y=(w-r_t)yv^*+y(v^*-s_{t'})$, and
$\norm{azb}_2\le\norm a\,\norm z_2\,\norm b$ for $a,b,z\in M$,
\[
 \norm{x-y}_2\le\lVert w-r_t\rVert\,\norm y_2\,\lVert v^*\rVert
 +\norm y_2\,\lVert v-s_{t'}\rVert\le4\theta\norm y_2 .
\]
By the triangle inequality, $L^2$-contractivity of $\E$, and
\eqref{eq:tail-lower},
\begin{align*}
 \norm{\E(x)}_2
 &\ge\norm{\E(y)}_2-\norm{\E(x-y)}_2
 \ge\norm{\E(y)}_2-\norm{x-y}_2\\
 &\ge\sqrt{2\ell+t+t'}\,\norm y_2-4\theta\norm y_2
 =\bigl(\sqrt{2\ell+t+t'}-4\theta\bigr)\norm x_2,
\end{align*}
with a positive coefficient by \eqref{eq:bufferconditions}.

\emph{Step 4: dimension count.} \Cref{lem:rectangle}, applied to
$\tilde r_t\le e_t$ and $\tilde s_{t'}\le f_{t'}$, and trace preservation under
partial isometries give
\begin{equation}\label{eq:distribution}
 \tau(r_t)\tau(s_{t'})\le d^{-1}\tau(e_t)\tau(f_{t'})
 \qquad(t,t'\ge0).
\end{equation}

\emph{Step 5: integration.} We first express the four $2$-norms in
\eqref{eq:tails} through the traces in \eqref{eq:distribution}. Since
$p$ commutes with $b$, the projections $p\one_B(b)$, for Borel sets
$B\subset\R$, define a finite positive Borel measure
$\mu_p(B)=\tau(p\one_B(b))$, and $\tau(pg(b))=\int g\,d\mu_p$ for every
bounded Borel function $g$. For $t\ge0$ and $\lambda\in\R$ we have
$(\lambda-c-\ell)_+>t$ if and only if $\lambda>c+\ell+t$. Hence
\eqref{eq:layer-cake}, applied to $(\lambda-c-\ell)_+$ and integrated
against $\mu_p$ (Tonelli), gives
\[
 \tau\bigl(p(b-c-\ell)_+^2\bigr)
 =\int_0^\infty2t\,\mu_p\bigl((c+\ell+t,\infty)\bigr)\,dt
 =\int_0^\infty2t\,\tau(r_t)\,dt .
\]
Since $(b-c-\ell)_+$ commutes with $p$, the left side is
$\norm{(b-c-\ell)_+p}_2^2$. In the same way, using
$\mu_q(B)=\tau(q\one_B(b))$ for $s_t$, and the spectral measure
$B\mapsto\tau(\one_B(a))$ of $a$ for $e_t$ and $f_t$,
\begin{align*}
 \norm{(c-\ell-b)_+q}_2^2&=\int_0^\infty2t\,\tau(s_t)\,dt,\\
 \norm{(a-c)_+}_2^2&=\int_0^\infty2t\,\tau(e_t)\,dt,\\
 \norm{(c-a)_+}_2^2&=\int_0^\infty2t\,\tau(f_t)\,dt .
\end{align*}
Now multiply \eqref{eq:distribution} by $4tt'=(2t)(2t')$ and integrate
over $(t,t')\in[0,\infty)^2$. By Tonelli, each side factors as a
product of two of the integrals above, so
\[
 \norm{(b-c-\ell)_+p}_2^2\,\norm{(c-\ell-b)_+q}_2^2
 \le d^{-1}\norm{(a-c)_+}_2^2\,\norm{(c-a)_+}_2^2 ,
\]
and taking square roots proves \eqref{eq:tails}.
\end{proof}

\subsection{The balanced threshold}
We will compare deviations above and below a threshold $c$. With $c=\alpha$
we have $\tau((a-\alpha)_+)=\tau((a-\alpha)_-)$, but
$\norm{(a-\alpha)_+}_2$ and $\norm{(a-\alpha)_-}_2$ need not agree, and
the argument would only control their product, which can be small while
$F$ is not. We therefore choose $c$ so that the two $2$-norms are equal.

\begin{lemma}[Balanced threshold]\label{lem:balance}
Let $p\in\Palpha$, $a=\E(p)$, $q=1-p$, and suppose
$F=\norm{a-\alpha1}_2^2>0$. There is $c\in[0,1]$ with
\[
 H:=\norm{(a-c)_+}_2^2=\norm{(c-a)_+}_2^2,
\]
and for any such $c$,
\begin{equation}\label{eq:balance}
 2H=F+(c-\alpha)^2\ge F>0,
\end{equation}
\begin{equation}\label{eq:covariance}
 \norm{(a-c)_+p}_2\ge\sqrt{\alpha H},\qquad
 \norm{(c-a)_+q}_2\ge\sqrt{(1-\alpha)H}.
\end{equation}
\end{lemma}
\begin{proof}
Since $0\le a\le1$, the continuous function
$c\mapsto\tau((a-c)_+^2)-\tau((c-a)_+^2)$ equals $\tau(a^2)\ge0$ at
$c=0$ and $-\tau((1-a)^2)\le0$ at $c=1$; choose a zero.

For \eqref{eq:balance}, note that $(a-c)_+(c-a)_+=0$ and
$a-c=(a-c)_+-(c-a)_+$, so $(a-c)^2=(a-c)_+^2+(c-a)_+^2$. Moreover
$\tau(a-\alpha1)=0$, since $\tau(a)=\tau(p)=\alpha$. Writing
$a-c=(a-\alpha1)+(\alpha-c)1$ and expanding,
\begin{align*}
 2H&=\tau\bigl((a-c)_+^2\bigr)+\tau\bigl((c-a)_+^2\bigr)=\tau\bigl((a-c)^2\bigr)\\
 &=\tau\bigl((a-\alpha1)^2\bigr)+2(\alpha-c)\,\tau(a-\alpha1)+(\alpha-c)^2
 =F+(c-\alpha)^2 .
\end{align*}

For \eqref{eq:covariance} we use Chebyshev's covariance inequality: if
$g$ is a non-decreasing function on $[0,1]$, then
$\tau(ag(a))\ge\alpha\,\tau(g(a))$. Indeed, since $\tau(a-\alpha1)=0$,
subtracting the scalar $g(\alpha)$ changes nothing, and
\[
 \tau(ag(a))-\alpha\,\tau(g(a))
 =\tau\bigl((a-\alpha1)(g(a)-g(\alpha)1)\bigr).
\]
The operator on the right is $\chi(a)$ with
$\chi(\lambda)=(\lambda-\alpha)(g(\lambda)-g(\alpha))$ for
$\lambda\in[0,1]$; since $g$ is non-decreasing, the two factors of
$\chi(\lambda)$ have the same sign, so $\chi\ge0$ and the trace is
nonnegative. Now take the non-decreasing function
$g(\lambda)=(\lambda-c)_+^2$. Since $g(a)\in N$ and $\E(p)=a$,
\[
 \norm{(a-c)_+p}_2^2=\tau(pg(a))=\tau\bigl(\E(p)g(a)\bigr)=\tau(ag(a))
 \ge\alpha\,\tau(g(a))=\alpha H.
\]
For the non-increasing function $g(\lambda)=(c-\lambda)_+^2$ the covariance
inequality reverses, so $\tau(qg(a))=\tau((1-a)g(a))\ge(1-\alpha)\tau(g(a))
=(1-\alpha)H$.
\end{proof}

The inequality \eqref{eq:covariance} is the density statement: $p$
carries at least an $\alpha$-fraction of $\norm{(a-c)_+}_2^2$, and $1-p$
at least a $(1-\alpha)$-fraction of $\norm{(c-a)_+}_2^2$.

\subsection{Conclusion of the proof}
The errors in the final comparison are controlled by the following
standard fact.

\begin{lemma}\label{lem:positive-part}
For self-adjoint $x,y\in M$, $\norm{x_+-y_+}_2\le\norm{x-y}_2$.
\end{lemma}
\begin{proof}
In the real Hilbert space $L^2(M)_{\sa}$, $x_+$ is the metric projection
of $x$ onto the closed convex cone of positive elements: for $z\ge0$,
$\langle x-x_+,z-x_+\rangle=-\tau(x_-z)\le0$. Metric projections onto
closed convex sets are nonexpansive.
\end{proof}

\begin{theorem}\label{thm:open-interval}
Let $N\subset M$ be an extremal inclusion of $\II$ factors with finite
index $d$, and let $\alpha$ satisfy $d\alpha(1-\alpha)>1$. Then
$\alpha\in\Lam(M,N)$.
\end{theorem}
\begin{proof}
Since $\alpha(1-\alpha)>d^{-1}$, we may choose $0<\delta<1$ with
$2\delta<\min(\sqrt\alpha,\sqrt{1-\alpha})$ and
\begin{equation}\label{eq:delta}
 (\sqrt\alpha-2\delta)(\sqrt{1-\alpha}-2\delta)>d^{-1/2}.
\end{equation}
Put $\eta=\delta^2/16$, and let $p\in\Palpha$ and $h$ be given by
\cref{cor:BP}. If $F=F(p)=0$, then $\E(p)=\alpha1$ and we are done.

Suppose $F>0$. As at the start of this section, $b=a-4F^{3/4}h$ commutes
with $p$, satisfies \eqref{eq:hessian}, and
$\kappa=\norm{a-b}\le4\eta F^{3/4}=\frac{\delta^2}4F^{3/4}$. Put
$\ell=\delta\sqrt F$. By \eqref{eq:F-bound}, $F\le1$, so
\[
 \frac\kappa\ell\le\frac\delta4F^{1/4}\le\frac14,\qquad
 \frac{4\kappa}\ell\le\delta F^{1/4}<\sqrt{2\delta}\,F^{1/4}=\sqrt{2\ell},
\]
and \eqref{eq:bufferconditions} holds. Let $c,H$ be as in
\cref{lem:balance}, and put
\[
 S_+=\norm{(b-c-\ell)_+p}_2,\qquad S_-=\norm{(c-\ell-b)_+q}_2 .
\]
\emph{Index side.} By \cref{lem:tails} and the choice of $c$,
\begin{equation}\label{eq:upperfinal}
 S_+S_-\le d^{-1/2}H .
\end{equation}
\emph{Density side.} By \cref{lem:positive-part},
$\norm{\bigl((b-c-\ell)_+-(a-c)_+\bigr)p}_2\le\norm{b-a-\ell}_2\le\kappa+\ell$,
so \eqref{eq:covariance} gives $S_+\ge\sqrt{\alpha H}-(\kappa+\ell)$, and
likewise $S_-\ge\sqrt{(1-\alpha)H}-(\kappa+\ell)$. Since $F\le1$,
$\delta<1$ and $F\le2H$,
\[
 \kappa+\ell\le\tfrac54\delta\sqrt F\le\tfrac{5\sqrt2}{4}\delta\sqrt H<2\delta\sqrt H,
\]
and therefore
\begin{equation}\label{eq:lowerfinal}
 S_+S_-\ge(\sqrt\alpha-2\delta)(\sqrt{1-\alpha}-2\delta)H .
\end{equation}
Since $H>0$, \eqref{eq:upperfinal} and \eqref{eq:lowerfinal} contradict
\eqref{eq:delta}. Hence $F=0$.
\end{proof}

\begin{proof}[Proof of \cref{thm:A}]
The condition $d\alpha(1-\alpha)>1$ is equivalent to
$\alpha\in(t_d,1-t_d)$, so \cref{thm:open-interval} gives
$(t_d,1-t_d)\subseteq\Lam(M,N)$. The endpoints $t_d$ and $1-t_d$, and
the description of $\Lam(M,N)$ outside $(t_d,1-t_d)$, are given by
\cref{thm:popa}(2), which holds for every inclusion of index $d\ge4$.
\end{proof}

This also proves \cref{cor:B}: part~(1) was observed after
\cref{thm:popa}, and part~(2) is \cref{thm:A}.

\begin{remark}
The strict inequality $d\alpha(1-\alpha)>1$ is needed for the margin
\eqref{eq:delta}. At the endpoints, Popa's argument is exact and of a
different nature \cite[proof of Thm.~5.2]{Popa89}. The tunnel of Jones
projections generates a pair of hyperfinite factors
$R_{-1}\subset R_0\subset M$ of
index $d$ that forms a commuting square with $N\subset M$
\cite[Prop.~1.8]{Popa89}. For $d>4$ this pair is a locally trivial
subfactor in the sense of \cref{def:locally-trivial}
\cite[Cor.~3.4]{Popa89} (see also \cite[Cor.~5.4]{PP86}), so its
defining projection $f$ satisfies
$E_{R_{-1}}(f)=t_d1$, and the commuting square gives
$\E(f)=t_d1$ \cite[Prop.~1.7]{Popa89}.
\end{remark}

\section{Popa's questions}\label{sec:final}

To apply \cref{cor:B} to Popa's questions we need that extremality is
preserved by compression and by the basic construction.

\begin{lemma}\label{lem:extremal-stable}
Let $N\subset M$ be an extremal inclusion of $\II$ factors with finite
index $d$.
\begin{enumerate}
\item For every $0\ne p\in\Proj(N)$, the inclusion $pNp\subset pMp$ is
 extremal of index $d$.
\item The inclusion $M\subset M_1=\basic MN$ is extremal of index $d$.
\end{enumerate}
\end{lemma}
The indices are \eqref{eq:jones-compression} and \cite{Jones83}.
Extremality in part~(1) is noted in the proof of
\cite[1.2.5(iv)]{Popa94}, and in part~(2) it is due to Pimsner and Popa
\cite{PP88}; see \cite[1.2.5(iii)]{Popa94}.

\begin{corollary}\label{cor:popa-problems}
Let $N\subset M$ be an extremal inclusion of $\II$ factors of finite
index $d$.
\begin{enumerate}
\item $\Lam(pMp,pNp)=\Lam(M,N)$ for every $0\ne p\in\Proj(N)$
 \textup{(}\cite[Problem~5.4.1]{Popa89}\textup{)}.
\item $\Lam(M_1,M)=\Lam(M,N)$
 \textup{(}\cite[Problem~5.4.2]{Popa89}\textup{)}.
\item $\Lam(M,N)$ is closed. It is uncountable if and only if $d>4$, and
 then it contains the interval $[t_d,1-t_d]$
 \textup{(}\cite[Problems~5.4.3 and~5.4.6]{Popa89}\textup{)}.
\end{enumerate}
\end{corollary}
\begin{proof}
By \cref{lem:extremal-stable}, $pNp\subset pMp$ and $M\subset M_1$ are
again extremal of index $d$, and by \cref{cor:B} the relative dimension
set of an extremal inclusion of finite index depends only on the index.
Part~(3) can be read off from the explicit sets in \cref{thm:popa} and
\cref{thm:A}, since $a_k\to t_d$.
\end{proof}
For $d\le4$, parts (1) and (2) hold without extremality, by
\cref{cor:B}(1). For $d>4$, part~(2) can fail without extremality
(\cref{rem:problem-542}); we do not know whether part~(1) can.

\begin{remark}[Unitary conjugacy]\label{rem:popa-questions}
Popa's question on unitary conjugacy remains open in the range of
\cref{thm:A}: if $N'\cap M=\C1$, $t_d<\alpha<1-t_d$, and
$f_1,f_2\in\Proj(M)$ satisfy $\E(f_1)=\E(f_2)=\alpha1$, are $f_1$ and $f_2$
conjugate by a unitary of $N$ \cite[Problem~5.4.5]{Popa89}? Outside
$[t_d,1-t_d]$ Popa proves that they are \cite[Thms.~5.1--5.2]{Popa89};
the endpoints are \cite[Problem~5.4.4]{Popa89}.
\end{remark}

\begin{remark}[Orthonormal unitaries]\label{rem:orthonormal-unitaries}
Elements $x_1,\dots,x_n\in M$ are \emph{$N$-orthonormal} if
$\E(x_j^*x_i)=\delta_{ij}1$. If $N\subset M$ is extremal with $d\ge4$,
\cref{thm:A} and \cref{thm:popa} (for $d=4$, where $t_4=1/2$) give
$p\in\Proj(M)$ with $\E(p)=\tfrac121$, so $u=2p-1$ is a self-adjoint
unitary with $\E(u)=0$, and $1,u$ are $N$-orthonormal. For inclusions
with principal graph $A_\infty$, Popa noted that finding even one
unitary $u\in M$ with $\E(u)=0$ was open \cite[\S3.5]{Popa19}, so
\cref{thm:A} answers this question. Popa also asks more generally how
many unitaries an orthonormal basis of $M$ over $N$ can contain.
\end{remark}

\begin{remark}[Partitions of unity]\label{rem:partitions}
Popa \cite[\S3.4]{Popa19} also considers the sets $\Lam_{k-1}(M,N)$ of
tuples $(\alpha_1,\dots,\alpha_k)$ for which there is a partition of
unity $p_1+\dots+p_k=1$ by projections with $\E(p_i)=\alpha_i1$; thus
$\Lam_1(M,N)$ corresponds to $\Lam(M,N)$. Are these sets determined by
the index for extremal inclusions of finite index? The proof of
\cref{thm:A} handles one projection at a time and does not produce
orthogonal families.
\end{remark}

\section{Non-extremal inclusions}\label{sec:nonextremal}

We now consider inclusions of finite index that need not be extremal.
For $d>4$, \cref{cor:C} shows that $\Lam(M,N)$ misses the middle interval
exactly for locally trivial inclusions. The remaining non-extremal
inclusions show both possible behaviors: $\Lam(M,N)$ can contain the
whole middle interval (\cref{ex:tensor}) or miss part of it
(\cref{ex:gap}).

\subsection{Locally trivial inclusions}
\begin{definition}[{\cite[Def.~3.1]{Popa89}}]\label{def:locally-trivial}
An inclusion $N\subset M$ of $\II$ factors is \emph{locally trivial} if
there are a projection $f\in M$ with $0<\tau(f)<1/2$ and an isomorphism
$\theta\colon fMf\to(1-f)M(1-f)$ such that
$N=\{x+\theta(x):x\in fMf\}$.
\end{definition}

\begin{remark}[Extremality cannot be dropped]\label{rem:locally-trivial}
If $M$ is hyperfinite, then every projection $f$ with $0<\tau(f)<1/2$
admits such an isomorphism $\theta$. Let $N\subset M$ be locally
trivial, with $f$ and $\theta$ as in \cref{def:locally-trivial}, and
put $s=\tau(f)$. Then $N\subset M$ has index
$d=s^{-1}+(1-s)^{-1}>4$ \cite[Exercise~9.18]{AP} and
$\Lam(M,N)=\{a_k,1-a_k:k\ge0\}\cup\{t_d,1-t_d\}$, so the middle interval
is missing \cite[\S3 and Thm.~5.3]{Popa89}. Here $f\in N'\cap M$ and
$Nf=fMf$, so $[fMf:Nf]=1\ne d\,\tau(f)^2=s/(1-s)$, and by
\eqref{eq:local-index} the inclusion is not extremal.
\end{remark}

For reducible inclusions, the proof of \cref{cor:C} is elementary; for
irreducible ones, it is \cref{thm:A}.

\begin{proof}[Proof of \cref{cor:C}]
If $N\subset M$ is locally trivial, the intersection is empty by
\cref{rem:locally-trivial}. Conversely, suppose it is empty. Then
$N'\cap M\ne\C1$, since otherwise $N\subset M$ is extremal and
\cref{thm:A} applies. Let $p\in N'\cap M$ be a projection with
$p\ne0,1$, and put $q=1-p$ and $s=\tau(p)$. Then $s\in\Lam(M,N)$
\cite[Prop.~1.4]{Popa89}. Applying \eqref{eq:local-index} to $p$ and
$q$ and using $\tau'(p)+\tau'(q)=1$ gives
\[
 d=\frac{[pMp:Np]}{s}+\frac{[qMq:Nq]}{1-s}\ge\frac1{s(1-s)} .
\]
If one of these local indices exceeded $1$, then $s(1-s)>d^{-1}$, that
is $s\in(t_d,1-t_d)$, contrary to assumption. So $pMp=Np$ and
$qMq=Nq$. Since $N$ is a factor, $x\mapsto xp$ and $x\mapsto xq$ are
isomorphisms of $N$ onto $pMp$ and $qMq$, so $\theta(xp)=xq$ defines an
isomorphism $\theta\colon pMp\to qMq$, and
$N=\{y+\theta(y):y\in pMp\}$. Moreover $d=1/(s(1-s))>4$ gives
$s\ne1/2$. Interchanging $p$ and $q$ if necessary, we get $s<1/2$, so
$N\subset M$ is locally trivial.
\end{proof}

\begin{remark}\label{rem:reducible}
The proof shows that if $N\subset M$ is not locally trivial, then every
projection of $N'\cap M$ other than $0,1$ has trace in $(t_d,1-t_d)$.
This covers every inclusion with $\dim(N'\cap M)\ge3$, since a locally
trivial inclusion has $N'\cap M=\C f+\C(1-f)$ \cite[Cor.~5.2.3(ii)]{Popa94}.

For a non-extremal inclusion that is not locally trivial, $\Lam(M,N)$
therefore meets the middle interval, but it need not contain it:
compare \cref{ex:tensor,ex:gap}. The description of $\Lam(M,N)$ for
non-extremal inclusions of finite index remains open. The proof of
\cref{thm:A} uses
extremality only through the existence of a tracial state on
$\mathcal S_0$ (\cref{thm:SEtrace}). In particular, $\mathcal S_0$ has
no tracial state when $N\subset M$ is locally trivial, since otherwise
that proof would give $(t_d,1-t_d)\subseteq\Lam(M,N)$.
\end{remark}

\subsection{An inclusion with the whole middle interval}
Our first example rests on two general observations. The first describes
how values of the relative dimension set combine under tensor products.

\begin{lemma}\label{lem:tensor}
Let $A\subset B$ and $P\subset Q$ be inclusions of $\II$ factors, and put
$N=A\mathbin{\overline\otimes}P$ and $M=B\mathbin{\overline\otimes}Q$. If
$\alpha\in\Lam(B,A)$ and $\beta,\beta'\in\Lam(Q,P)$, then
$\alpha\beta+(1-\alpha)\beta'\in\Lam(M,N)$.
\end{lemma}
\begin{proof}
Let $p\in B$ and $p_1,p_2\in Q$ be projections with $E_A(p)=\alpha1$,
$E_P(p_1)=\beta1$ and $E_P(p_2)=\beta'1$. Then
$p\otimes p_1+(1-p)\otimes p_2$ is a projection of $M$, and since
$\E=E_A\otimes E_P$, its expectation is
$(\alpha\beta+(1-\alpha)\beta')1$.
\end{proof}

The second observation concerns inclusions that are isomorphic to their
amplifications. Let $R$ be the hyperfinite $\II$ factor. In Popa's
terminology, $N\subset M$ \emph{splits} $R$ if
$N=N_0\mathbin{\overline\otimes}R$ and
$M=M_0\mathbin{\overline\otimes}R$ for an inclusion $N_0\subset M_0$ of
$\II$ factors. For such inclusions, Popa obtains from one value of
$\Lam(M,N)$ in $(t_d,1/2)$ a sequence of values decreasing to $t_d$
\cite[Prop.~5.5]{Popa89}. Applied to a whole interval, and to values
above $1/2$ as well, the same argument fills the middle interval.

\begin{lemma}\label{lem:propagate}
Let $N\subset M$ be an inclusion of $\II$ factors with
$4<d=[M:N]<\infty$ that splits $R$. If
$[2/d,1-2/d]\subseteq\Lam(M,N)$, then
$[t_d,1-t_d]\subseteq\Lam(M,N)$.
\end{lemma}
\begin{proof}
Put $S(x)=1/(d(1-x))$, the recursion in \eqref{eq:an}, and $T=S\circ S$.
The proof of \cite[Prop.~5.5]{Popa89} shows that $T(x)\in\Lam(M,N)$ for
every $x\in\Lam(M,N)\cap(t_d,1-t_d)$: it applies
\cite[Prop.~4.1]{Popa89} twice and uses that, in the Jones tower
$N\subset M\subset M_1\subset M_2$, the inclusion $M_1\subset M_2$ is an
amplification of $N\subset M$ \cite[1.3.3]{Popa94}, and the assumption
$\alpha<1/2$ made there is not needed for this.

On $(t_d,1-t_d)$ we have $dx(1-x)>1$, that is $S(x)<x$. Since $S$ is
increasing with $S(t_d)=t_d$, the map $T$ is continuous and increasing
there, with $t_d<T(x)<x$, and a direct computation gives
$T(1-2/d)=2/d$. Since $d>4$, we have $2/d\in(t_d,1/2)$.
Let $b_0=2/d$ and $b_{n+1}=T(b_n)$, so that $b_n$ decreases to $t_d$.
Then $T$ maps $[b_0,1-b_0]$ onto $[b_1,b_0]$, and $[b_n,b_{n-1}]$ onto
$[b_{n+1},b_n]$ for $n\ge1$, so all these intervals lie in $\Lam(M,N)$.
Hence $(t_d,1-2/d]\subseteq\Lam(M,N)$. The invariance of $\Lam(M,N)$
under $\alpha\mapsto1-\alpha$ gives $[2/d,1-t_d)\subseteq\Lam(M,N)$, and
$t_d,1-t_d\in\Lam(M,N)$ by \cref{thm:popa}.
\end{proof}

\begin{example}\label{ex:tensor}
The conclusion of \cref{thm:A} can hold without extremality, even when
$N'\cap M$ is two-dimensional, so that its projections contribute only
two values to the middle interval. Let $A\subset B$ be a
locally trivial inclusion of $\II$ factors with $B$ hyperfinite and
$\tau(f)=1/3$, so that $[B:A]=9/2$ and $A'\cap B=\C f+\C(1-f)$
(\cref{rem:locally-trivial,rem:reducible}). Let $P_0\subset Q_0$ be an
irreducible inclusion of $\II$ factors with $[Q_0:P_0]=5$, for instance
$Q_0=P_0\rtimes\Z/5\Z$ for an outer action, and put
$P=P_0\mathbin{\overline\otimes}R$ and
$Q=Q_0\mathbin{\overline\otimes}R$, so that $P\subset Q$ is irreducible
of index $5$. Put
$N=A\mathbin{\overline\otimes}P$, $M=B\mathbin{\overline\otimes}Q$ and
$r=f\otimes1$. Then $d=[M:N]=45/2$ and
$N'\cap M=(A'\cap B)\mathbin{\overline\otimes}(P'\cap Q)=\C r+\C(1-r)$,
so the traces of its projections other than $0,1$ are $1/3$ and $2/3$.
Since $Af=fBf$, the
inclusion $Nr\subset rMr$ is
$fBf\mathbin{\overline\otimes}P\subset fBf\mathbin{\overline\otimes}Q$,
so $[rMr:Nr]=5$, and likewise $[(1-r)M(1-r):N(1-r)]=5$. Both local
indices exceed $1$, so $N\subset M$ is not locally trivial, and since
$5\ne d\,\tau(r)^2=5/2$ it is not extremal.

By \cref{thm:A}, $[t_5,1-t_5]\subseteq\Lam(Q,P)$, and by
\cref{thm:popa}, $\Lam(B,A)$ contains $1/3=t_{9/2}$, as well as $2/9$
and $7/9$, the values $a_1$ and $1-a_1$ for index $9/2$; this holds for
every inclusion of index $9/2$.
\Cref{lem:tensor} with $\beta'=0$ gives
$\alpha[t_5,1-t_5]\subseteq\Lam(M,N)$ for these three values of
$\alpha$. These intervals overlap and cover
$[2t_5/9,7(1-t_5)/9]$, which contains $1/2$. Since $\Lam(M,N)$ is
invariant under $\alpha\mapsto1-\alpha$,
\[
 [2t_5/9,1-2t_5/9]\subseteq\Lam(M,N).
\]
Now $N\subset M$ splits $R$, with
$N_0=A\mathbin{\overline\otimes}P_0$ and
$M_0=B\mathbin{\overline\otimes}Q_0$, and $2t_5/9\approx0.061$ is
smaller than $2/d\approx0.089$. So \cref{lem:propagate} gives
$[t_d,1-t_d]\subseteq\Lam(M,N)$. By \cref{thm:popa}, $\Lam(M,N)$ is
therefore the set in \cref{thm:A}, although $N\subset M$ is not
extremal.
\end{example}

\subsection{An inclusion with a gap in the middle interval}
The opposite behavior comes from projections of $N'\cap M$ of local
index $1$.

\begin{lemma}\label{lem:gap}
Let $N\subset M$ be an inclusion of $\II$ factors, and let
$r\in N'\cap M$ be a projection with $Nr=rMr$. If $s=\tau(r)$, then
$\Lam(M,N)\cap(1-s,s)=\emptyset$.
\end{lemma}
\begin{proof}
Let $p\in M$ be a projection with $\E(p)=\alpha1$ and $\alpha>1-s$; we
show that $\alpha\ge s$. By Kaplansky's formula,
$\tau(p\wedge r)\ge\tau(p)+\tau(r)-1=\alpha+s-1>0$, so
$g=p\wedge r\ne0$. Since $g\in rMr=Nr$ and $x\mapsto xr$ is an
isomorphism of $N$ onto $Nr$, we have $g=zr$ for a nonzero projection
$z\in N$. As $\E(r)=s1$, this gives $\E(g)=z\E(r)=sz$. Since $g\le p$,
we get $sz=\E(g)\le\E(p)=\alpha1$, and so $s\le\alpha$.
\end{proof}

For a locally trivial inclusion, taking $r=1-f$ gives $s=1-t_d$, which
recovers the fact that the middle interval is missing.

\begin{example}\label{ex:gap}
A non-extremal inclusion that is not locally trivial can miss part of
the middle interval. Let $M$ be hyperfinite, let $e\in\Proj(M)$ with
$\tau(e)=2/3$, and put $f=1-e$. Let $P\subset fMf$ be an irreducible
subfactor of index $2$, choose an isomorphism $\theta\colon eMe\to P$,
and put $N=\{x+\theta(x):x\in eMe\}$, so that $Ne=eMe$ and $Nf=P$.

Then $N'\cap M=\C e+\C f$. Indeed, if $z\in N'\cap M$, then
$eze\in\C e$ and $fzf\in P'\cap fMf=\C f$. Moreover $v=ezf$ satisfies
$xv=v\theta(x)$ for $x\in eMe$, so $vv^*\in\C e$ and $v^*v\in\C f$;
since $\tau(e)\ne\tau(f)$, this forces $v=0$, and likewise $fze=0$. By
\eqref{eq:local-index}, as in the proof of \cref{cor:C},
$d=[M:N]=1/(2/3)+2/(1/3)=15/2$. Since $[fMf:Nf]=2$, the inclusion is not
locally trivial, and since $[eMe:Ne]=1\ne d\,\tau(e)^2=10/3$, it is not
extremal.

The values $1/3$ and $2/3$ lie in $\Lam(M,N)$ \cite[Prop.~1.4]{Popa89},
and inside the middle interval, as $t_d\approx0.158$. On the other hand,
\cref{lem:gap} with $r=e$ gives $\Lam(M,N)\cap(1/3,2/3)=\emptyset$.
Since a locally trivial inclusion of index $15/2$ misses the whole
middle interval, $\Lam(M,N)$ is not determined by the index for
non-extremal inclusions.
\end{example}

\begin{remark}[Popa's Problem 5.4.2]\label{rem:problem-542}
The inclusion of \cref{ex:gap} shows that \cref{cor:popa-problems}(2)
fails without extremality. Let $f'=JfJ\in M'\cap M_1$. The anti-isomorphism
$x\mapsto Jx^*J$ of $N'$ onto $M_1$ carries $M'f\subset fN'f$ onto
$Mf'\subset f'M_1f'$, so
$\tau_{M_1}(f')=\tau'(f)=[fMf:Nf]/(d\,\tau(f))=4/5$ by
\eqref{eq:local-index}, and
$[f'M_1f':Mf']=d\,\tau(f)\tau'(f)=2$ \cite[1.2.1(a)]{Popa94}. So there is
a projection $g\le f'$ whose expectation onto $Mf'$ is $\tfrac12f'$
\cite[1.2.2--1.2.3]{Popa94}, and comparing traces against elements of
$M$ gives $E_M(g)=\tfrac45\cdot\tfrac121=\tfrac251$. Thus
$2/5\in\Lam(M_1,M)$, while $2/5\notin\Lam(M,N)$. This answers the first
part of \cite[Problem~5.4.2]{Popa89} in the negative; its second part,
for $N'\cap M=\C1$, is answered positively by
\cref{cor:popa-problems}(2).
\end{remark}

\bibliographystyle{amsplain}
\bibliography{refs}

\end{document}